\documentstyle{article}
\title{\bf Geometrical interpretation of solutions\\
of certain PDEs}
\author{C. Udri\c ste, M. Neagu}
\date{}

\begin{document}
\maketitle
\begin{abstract}
In \S1 the authors define the notion of harmonic map between
two ge\-ne\-ra\-lized Lagrange spaces. In \S2 it is proved that for certain systems of
differential or partial differential equations, the solutions belong to a class
of harmonic maps between two generalized Lagrange spaces.
\S3 describes the main properties of the generalized Lagrange spaces
constructed in \S2. These spaces, being con\-ve\-ni\-ent relativistic models, allow
us to write the Maxwell's and Einstein's equations.
\end{abstract}

\section{Introduction}

\hspace{5mm}Looking for generalizing a Poincar\'e problem, Sasaki tried to
find a Riemannian metric on a manifold
$M$ such that the orbits of an arbitrary vector field $X$ should be geodesics.
This attempt was a failure, but Sasaki discovered the well known
almost contact metric structures on a manifold of odd dimension \cite{8}.
After the introduction of generalized Lagrange structures
\cite{4}, the problem of Poincar\'e-Sasaki was reconsidered by the first
author \cite{9,10,11}. He succeded to discover a Lagrange structure on $M$,
depending of the given vector field $X$ and a $(1,1)$-tensor field
built using $X$, a metric $g$, and the covariant derivative induced by $g$,
such that the $C^2$ orbits should belong to a class of geodesics. Moreover,
replacing the system of ODEs of the orbits of $X$ by a system of PDEs and the
notion of geodesic by the notion of harmonic map, same open general problems
appear \cite{12}, namely

1) There exist Lagrange type structures such that the solutions of
certain PDEs of order one should be {\it harmonic maps}?

2) What is a {\it harmonic map} between two generalized Lagrange spaces?

Using the notion of {\it direction dependent harmonic map} between a Riemannian manifold
and a generalized Lagrange space, a partial answer to the Udri\c ste's
questions was offered by the second author \cite{6}.

Let us introduce, in a natural way, the notion of harmonic
map between two Lagrange spaces $(M,\;g_{\alpha\beta}(a,b))$ and
$(N,\;h_{ij}(x,y))$, where $M$ (resp. $N$) has the dimensions $m$ (resp. $n$)
and $(a,b)=(a^\mu ,b^\mu )$ (resp. $(x,y)=(x^k,y^k))$ are coordinates on $TM$
(resp. $TN$).

{\bf Definition.} On $M\times N$, a tensor $P$ of type $(1,2)$
with all components null excepting $P^\beta_{\alpha i}(a, x)$ and
$P^j_{\alpha i}(a, x)$, where $\alpha ,\beta =\overline{1,m},\;i,j=\overline{1,n}$,
is called {\it tensor of connection}.

Assume that the manifold $M$ is connected, compact, orientable and endowed also
with a Riemannian metric $\varphi_{\alpha\beta}$. This fact ensures the existence
of a volume element on $M$. In these conditions, we can define the
$\left(
\begin{array}{ccc}
&P&\\
g&\varphi&h
\end{array}\right)$
{\it -energy functional},
$
E\stackrel{\hbox{not}}=E^P_{g\varphi h}:C^\infty(M,N)\to R
$,
$$
E^P_{g\varphi h}(f)=\displaystyle{{1\over 2}\int_Mg^{\alpha\beta}(a,b)
h_{ij}(f(a),y)f^i_\alpha f^j_\beta\sqrt\varphi da,}
$$
where
$\left\{\begin{array}{l}\medskip
\displaystyle{
f^i=x^i(f),\;f^i_\alpha ={\partial f^i\over\partial a^\alpha},\;
\varphi=det(\varphi_{\alpha\beta})},\\\medskip
\displaystyle{b(a)=b^\gamma (a)
\left.{\partial\over\partial a^\gamma}\right\vert_a
\stackrel{\hbox{def}}=\varphi^{\alpha\beta}(a)f^i_\alpha (a)P^\gamma_{\beta i}(a,f(a))
\left.{\partial\over\partial a^\gamma}\right\vert_a
},\\
\displaystyle{
y(f(a))=y^k(a)\left.{\partial\over\partial x^k}\right\vert_{f(a)}
\stackrel{\hbox{def}}=\varphi^{\alpha\beta}(a)f^i_\alpha (a)P^k_{\beta i}(a,f(a))
\left.{\partial\over\partial x^k}\right\vert_{f(a)}}.
\end{array}\right.$

{\bf Definition.} A map $f\in C^\infty(M,N)$ is called
$\left(
\begin{array}{ccc}
&P&\\
g&\varphi&h
\end{array}\right)$
{\it -harmonic} if $f$ is a critical point for the functional
$E^P_{g\varphi h}$.

The naturalness of the preceding definitions comes from the following particular
cases:

i) If $g_{\alpha\beta}(a,b)=\varphi_{\alpha\beta}(a)$ and $h_{ij}(x,y)=h_{ij}(x)$
are Riemannian metrics, it recovers the classical definition of a harmonic map
between two Riemannian manifolds \cite{2,3}.

ii) If $M=[a,b]\subset R,\;\varphi_{11}=g_{11}=1$ and $P=(P^1_{1i},\delta^k_i)$
we shall find \linebreak $C^\infty (M,N)=\{c:[a,b]\to N\vert\;c-C^\infty \hbox{differentiable}
\}\stackrel{\hbox{not}}=\Omega_{a,b}(N)$ and the energy functional will be
$$
E^P_{11h}(c)=\displaystyle{{1\over 2}\int_a^bh_{ij}(c(t),\dot c(t)){dc^i\over
dt}{dc^j\over dt}dt,\;\forall c\in\Omega_{a,b}(N).}
$$
In conclusion, the
$\left(
\begin{array}{ccc}
&P&\\
1&1&h
\end{array}\right)$-harmonic curves are exactly the geodesics of the generalized
Lagrange space $(N,h_{ij}(x,y))$ \cite{4}.

iii) If we take $N=R,\;h_{11}=1$ and $P=(\delta^\alpha_\beta,P^1_{\beta 1})$ we
shall obtain
$C^\infty(M,N)=\nolinebreak{\cal F}(M)$ and the energy functional becomes
$$
E^P_{g\varphi 1}(f)=\displaystyle{
{1\over 2}\int_Mg^{\alpha\beta}(a,grad_\varphi f)f_\alpha f_\beta\sqrt\varphi
da,\;\forall f\in{\cal F}(M).}
$$

Obviously, the Euler-Lagrange equations of the energy functional $E$ will
be the equations of {\it harmonic maps}, that is,
$$\displaystyle{\sqrt\varphi{\partial L\over\partial f^i}+{\partial\over
\partial a^\alpha}\left(\sqrt\varphi{\partial L\over\partial f^i_\alpha}
\right)=0,\;\forall i=\overline{1,n},}\leqno(H)
$$
where $\displaystyle{L(a^\alpha,f^i,f^i_\alpha)={1\over 2}g^{\gamma\mu}(a^\nu,
b^\nu)h_{kl}(x^p,y^p)f_\gamma^kf^l_\mu}$.

In the particular cases when the metric tensors are
$g_{\alpha\beta}(a,b)=e^{-2\sigma(a,b)}\varphi_{
\alpha\beta}(a)$ and $h_{ij}(x,y)=e^{2\tau(x,y)}\psi_{ij}(x)$, where
$\sigma:TN\to R,\;\tau:TN\to R$ are smooth functions and $\psi_{ij}$ is a
pseudo-Riemannian metric on $N$, we shall obtain
$$\left\{\begin{array}{l}
\displaystyle{{\partial L\over\partial x^i}=e^{2\sigma+2\tau}\varphi^{\gamma
\mu}\varphi^{\delta\varepsilon}\psi_{kl}\left[{\partial P^\nu_{\varepsilon p}
\over\partial x^i}{\partial\sigma\over\partial b^\nu}+{\partial P^j_
{\varepsilon p}\over\partial x^i}{\partial\tau\over\partial y^j}\right]x^p_\delta
x^k_\gamma x^l_\mu+{1\over 2}g^{\gamma\mu}{\partial h_{kl}\over\partial x^i}x^k_
\gamma x^l_\mu}\\
\displaystyle{{\partial L\over\partial x^i_\alpha}=e^{2\sigma+2\tau}\left\{
\varphi^{\gamma\mu}\varphi^{\alpha\varepsilon}\psi_{kl}\left[P^\nu_{\varepsilon
i}{\partial\sigma\over\partial b^\nu}+P^j_{\varepsilon i}{\partial\tau\over
\partial y^j}\right]x^k_\gamma x^l_\mu+\varphi^{\gamma\alpha}\psi_{ik}x^k_
\gamma\right\}.}
\end{array}\right.
$$

These expressions will be reduced if we consider the following more particular
cases:

1) $\sigma=\sigma(a)$ and $P=(P^\alpha_{\beta i},A_\beta(a)\delta^i_j)$, where
$\{A_\beta\}$ are the components of a covector $A$ on $M$. In this situation, we
shall obtain
$$\left\{
\begin{array}{l}\medskip
\displaystyle{
{\partial L\over\partial x^i}={1\over 2}g^{\gamma\mu}{\partial
h_{kl}\over\partial x^i}x^k_\gamma x^l_\mu}\\
\displaystyle{
{\partial L\over\partial x^i_\alpha}=e^{2\sigma+2\tau}\left\{
\varphi^{\gamma\mu}\varphi^{\alpha\varepsilon}\psi_{kl}A_\varepsilon
{\partial\tau\over\partial y^i}x^k_\gamma x^l_\mu+\varphi^{\gamma\alpha}
\psi_{ik}x^k_\gamma\right\}.}
\end{array}\right.\leqno(*)
$$

2) $\tau=\tau(x)$ and $P=(\delta^\alpha_\beta\xi_i(x), P^\alpha_{\beta i})$,
where $\{\xi_i\}$ are the components of an 1-form $\xi$ on $N$. Now, we shall
find
$$\left\{
\begin{array}{l}\medskip
\displaystyle{
{\partial L\over\partial x^i}=e^{2\sigma+2\tau}
\varphi^{\gamma\mu}\varphi^{\delta\varepsilon}\psi_{kl}{\partial\xi_p\over
\partial x^i}{\partial\sigma\over\partial b^\varepsilon}x^p_\delta x^k_\gamma
x^l_\mu+{1\over 2}g^{\gamma\mu}{\partial h_{kl}\over\partial x^i}x^k_\gamma
x^l_\mu,}\\
\displaystyle{
{\partial L\over\partial x^i_\alpha}=e^{2\sigma+2\tau}\left\{
\varphi^{\gamma\mu}\varphi^{\alpha\varepsilon}\psi_{kl}
{\partial\sigma\over\partial b^\varepsilon}\xi_ix^k_\gamma x^l_\mu+\varphi^
{\gamma\alpha}\psi_{ik}x^k_\gamma\right\}.}
\end{array}\right.\leqno(**)
$$

\section{Geometrical interpretation}

\hspace{5mm}By the above notions, we shall offer some beautiful geometrical
interpretations of $C^2$ solutions of certain PDEs of order one.

We start with a smooth map $f\in C^\infty(M,N)$. This map induces the
following tensor
$\displaystyle{
\left.\delta f\stackrel{\hbox{not}}=f^i_\alpha da^\alpha
\otimes{\partial\over\partial y^i}\right\vert_{f(x)}\in \Gamma(T^*M\otimes
f^{-1}(TN))}$. On $M\times N$, let $T$ be a tensor of type $(1,1)$ with all
components null excepting $(T^i_{\alpha} )_{i=\overline{1,n}\\
\atop
\alpha=\overline{1,m}}$. These objects determine the system of PDEs,
$$\displaystyle{
\delta f=T\;\hbox{expressed locally by}\;{\partial f^i\over\partial a^\alpha}
=T^i_\alpha(a,f).
}\leqno(E)$$

If $(M,\varphi_{\alpha\beta})$ and $(N,\psi_{ij})$ are Riemannian
manifolds we can build a scalar product on $\Gamma(T^*M\otimes f^{-1}(TN))$,
namely
$<T,S>=\varphi^{\alpha\beta}\psi_{ij}T^i_\alpha S^j_\beta$, where
$\displaystyle{T=T^i_\alpha da^\alpha\otimes{\partial\over\partial y^i}}$ and
$\displaystyle{S=S^j_\beta da^\beta\otimes{\partial\over\partial y^j}}$.
Obviously, the Cauchy-Schwartz inequality
$$
<T,S>\leq\Vert S\Vert^2\Vert T\Vert^2,\; \forall S,T\in\Gamma(T^*M\times
f^{-1}(TN)),
$$
is an equality iff there exists ${\cal K}\in {\cal F}(M)$ such that
$T={\cal K}S$.

In these conditions, we prove the following

{\bf Theorem.} {\it
If $(M,\varphi),(N,\psi)$ are Riemannian manifolds and the smooth map
\linebreak $f\in C^\infty(M,N)$ is solution of the system $(E)$, then $f$ is an
extremal of functional
$$
{\cal L}_T:C^\infty(M,N)\backslash\{\exists a\in M\;\mbox{such that}\;
<\delta f,T>(a)=0\}\to R_+,
$$
$$
\displaystyle{{\cal L}_T(f)={1\over 2}\int_M{\Vert\delta f\Vert^2\Vert T\Vert^2
\over <\delta f,T>^2}\sqrt{\varphi}da={1\over 2}\int_M{\Vert T\Vert^2\over <
\delta f, T>^2}\varphi^{\alpha\beta}\psi_{ij}f^i_\alpha f^j_\beta\sqrt{\varphi}da}.
$$}

{\bf Proof.} Let $f$ be an arbitary map from the definition domain of ${\cal L}_T$.
Applying the above Cauchy-Schwarz inequality, we get
$$
\displaystyle{
{\cal L}_T(f)={1\over 2}\int_M{\Vert\delta f\Vert^2\Vert T\Vert^2\over
<\delta f,T>^2}\sqrt\varphi da\geq{1\over 2}\int_M\sqrt\varphi da={1\over 2}Vol_\varphi
(M).
}$$
Obviously, if $f$ is solution of the system $(E)$ it follows
$\displaystyle{{\cal L}_T(f)={1\over 2}Vol_\varphi(M)}$, that is, $f$ is a
global minimum point of the functional ${\cal L}_T$. In conclusion, the map
$f$ verifies the Euler-Lagrange equations of ${\cal L}_T\;\Box$.

Generally, the global minimum points of the functional ${\cal L}_T$ are
solutions of the system $\delta f={\cal K}T$, where ${\cal K}\in{\cal F}(M)$.
They are not necessarily solutions of initial system \nolinebreak $(E)$.

Now, we remark that, in certain particular cases, the functional ${\cal L}_T$ becomes
exactly a functional of type
$\left(
\begin{array}{ccc}
&P&\\
g&\varphi&h
\end{array}\right)$-energy and, consequently, the Euler-Lagrange equations
reduce to equations of harmonic maps.
This fact allows the following geometrical interpretations:

{\bf 1. Orbits}

Taking $M=([a,b],1)$ and $T=\xi\in\Gamma(c^{-1}(TN))$, the
PDEs system $(E)$ reduces to the system of orbits
$$\displaystyle{{dc^i\over dt}=\xi^i(c(t)),\;c:[a,b]\to N},$$
and the functional ${\cal L}_\xi$ comes to
$$\displaystyle{{\cal L}_\xi(c)={1\over 2}\int^b_a{\Vert\xi\Vert^2_\psi\over
[\xi^b(\dot c)]^2}\psi_{ij}{dc^i\over dt}{dc^j\over dt}dt},
$$
where $\xi^b=\xi_idx^i=\psi_{ij}\xi^jdx^i$.
Hence the functional ${\cal L}_\xi$ is a
$\left(
\begin{array}{ccc}
&P&\\
1&1&h
\end{array}\right)$-energy,
where the Lagrange metric tensor
$h_{ij}:TN\backslash\{y\vert \xi^b(y)=0\}\to R$
is defined by
$$\displaystyle{
h_{ij}(x,y)={\Vert \xi\Vert^2_\psi\over[\xi^b(y)]^2}\psi_{ij}(x)=
\psi_{ij}(x)\exp{\left[2\ln{\Vert \xi\Vert_\psi\over\vert\xi^b(y)\vert}\right]}
}.$$
This case was studied, in other way, by Udri\c ste \cite{9}-\cite{12}.

Replacing $\sigma=0,\;\tau(x,y)=\ln(\Vert\xi\Vert_\psi/
\vert\xi^b(y)\vert),\;\varphi_{11}=1$ and $A_1=1$
in the equations $(*)$,the following equations
will be the equations of these harmonic curves,
$$
\displaystyle{{\partial L\over\partial c^i}+{d\over dt}{\partial
L\over\partial \dot c^i}=0,\;\forall i=\overline{1,n},}
$$
where $\displaystyle{{\partial L\over \partial c^i}={1\over 2}{h_{kl}\over\partial
c^i}\dot c^k\dot c^l}$ and
$\displaystyle{{\partial L\over\partial\dot c^i}=e^{2\tau}\left\{
\psi_{kl}{\partial\tau\over\partial\dot c^i}\dot c^k\dot c^l+\psi_{ik}\dot c^k
\right\}.}$
\medskip

{\bf 2. Pfaff systems}

If we put $N=(R,1)$ and $T=A\in\Lambda^1(T^*M)$, the PDEs system $(E)$
becomes the Pfaffian system
$$df=A,\;f\in{\cal F}(M)$$
and the functional ${\cal L}_T$ is
$$
\displaystyle{{\cal L}_A(f)={1\over 2}\int_M{\Vert A\Vert^2_\varphi\over[A(
grad_\varphi f)]^2}\varphi^{\alpha\beta}f_\alpha f_\beta\sqrt{\varphi}da}.
$$
Consequently, the functional ${\cal L}_A$ is a
$\left(
\begin{array}{ccc}
&P&\\
g&\varphi&1
\end{array}\right)$-energy, where \linebreak
$g_{\alpha\beta}: TM\backslash\{b\vert A(b)=0\}\to R$ is defined by
$$\displaystyle{
g_{\alpha\beta}(a,b)={[A(b)]^2\over\Vert A\Vert^2_\varphi}\varphi_{\alpha\beta}(a)
=\varphi_{\alpha\beta}(a)
\exp{\left[2\ln{\vert A(b)\vert\over\Vert A\Vert_\varphi}\right]}.}$$

The form of harmonic maps equations are obtained, in this case, replacing
$\tau=0,\;\sigma(a,b)=\ln(\Vert A\Vert_\varphi/\vert A(b)\vert),\;\psi_{11}=
h_{11}=1$ and $\xi_1=1$ in $(**)$. These will be the equations
of harmonic maps $(H)$ with $n=1$, where
$$\displaystyle{{\partial L\over\partial f_\alpha}=e^{2\tau}\left\{
\varphi^{\gamma\mu}\varphi^{\alpha\varepsilon}{\partial\sigma\over\partial b^
\varepsilon}f_\gamma f_\mu+\varphi^{\gamma\alpha}f_\gamma\right\},
\;{\partial L\over\partial f}=0}.$$
\medskip

{\bf 3. Pseudolinear functions}

We suppose that $T^k_\beta (a,x)=\xi^{k}(x)A_\beta (a)$, where $\xi^k$
is vector field on $N$ and $A_\beta$ is 1-form on $M$. In this case the
functional ${\cal L}_T$ is expressed by
$${\cal L}_T(f)
\displaystyle{={1\over 2}\int_M{\Vert\xi\Vert^2_\varphi\Vert A
\Vert^2_\varphi\over[A(b)]}\varphi^{\alpha\beta}\psi_{ij}f^i_\alpha f^j_\beta
\sqrt\varphi da=}
$$
$$\displaystyle{={1\over 2}\int_Mg^{\alpha\beta}(a,b)h_{ij}(f(a))f^i_\alpha f^j
_\beta\sqrt\varphi da},
$$
where $P^\gamma_{i\beta}(x)=\delta^\gamma_\beta\xi_i(x),\;b^\gamma=
\varphi^{\alpha\beta}f^i_\alpha P^\gamma_{i\beta},\;h_{ij}(x)=\Vert
\xi\Vert^2_\psi\psi_{ij}(x)$ and the Lagrange metric tensor
$g_{\alpha\beta}:TM\backslash\{b\vert A(b)=0\}\to R$ is defined by
$$\displaystyle{g_{\alpha\beta}(a,b)={[A(b)]^2\over\Vert A\Vert^2_\varphi}
\varphi_{\alpha\beta}(a)=\varphi_{\alpha\beta}(a)
\exp{\left[2\ln{\vert A(b)\vert\over\Vert A\Vert_\varphi}\right]}
}.$$
It follows that the functional ${\cal L}_T$ becomes a
$\left(
\begin{array}{ccc}
&P&\\
g&\varphi&h
\end{array}\right)$-energy.

The equations of harmonic maps can be computed, putting 
$$
\tau=\ln\Vert\xi\Vert_\psi,\;\sigma(a,b)=\ln(\Vert A\Vert_\varphi/\vert
A(b)\vert),\;P^\gamma_{\beta i}=\delta^\gamma_\beta\xi^b_i(x),
$$
in $(**)$, where $\xi^b_i=\psi_{ij}\xi^j$.

In the particular case when we have
$M=(R^n, \varphi=\delta)$ and $N=(R,\psi=1)$, supposing that
$(grad\;f)(a)\neq 0,\;\forall a\in M$, the solutions
of the above system are the well known {\it pseudolinear functions} \cite{7}.
These functions have the property that all hypersurfaces of constant level
$M_{f(a)}$ are totally geodesic \cite{7}. Consequently, the pseudolinear
functions are examples of harmonic maps between the generalized Lagrange
spaces $(M,g_{\alpha\beta}(a,b)=\delta_{\alpha\beta}
\{[A(b)]^2/\Vert A\Vert^2\})$ and $(N,h(x)=\xi^2(x))$. For example, the
function $f(a)=e^{<v,a>+w}$, where $v\in M,\;w\in R$, is solution for above
system with $\xi(a)=1$ and $A(f(a))=f(a)v$.
\medskip

{\bf 4. Continuous groups of transformations}

The fundamental PDEs system of the group having the infinitesimal generators $\xi_r$
are
$$
{\partial f^i\over\partial a^\alpha}=\sum^t_{r=1}\xi_r^i(f)A^r_\alpha(a),
$$
where $\{\xi_r\}_{r=\overline{1,t}}\subset{\cal X}(N)$ are vector
fields on $N$ and $\{A^r\}_{r=\overline{1,t}}\subset\Lambda^1(M)$ is a familly
of covector fields on $M$. The geometrical interpretation of solutions via
harmonic maps theory is still an open problem, though the Lagrangian
$$
\displaystyle{L(a^\alpha,f^i,f^i_\alpha)={1\over 2}g^{\gamma\mu}(a^\nu,
b^\nu)h_{kl}(x^p,y^p)f_\gamma^kf^l_\mu}
$$
is what we need.

\section{Maxwell and Einstein equations}

\hspace{5mm}Finally, we remark that, in all above cases, the solutions of the system
$\delta f=T$ are harmonic maps between generalized Lagrange spaces of type
$(M^n,e^{2\sigma(x,y)}\gamma_{ij}(x))$, where $\sigma:TM\backslash\{\mbox
{Hyperplane}\}\to R$ is a smooth function. These spaces, endowed with the
non-linear connection $N^i_j(x,y)=\Gamma^i_{jk}(x)y^k$, where $\Gamma^i_{jk}(x)$
are the Christoffel symbols for the Riemannian metric $\gamma_{ij}(x)$,
verify a constructive
axiomatic formulation of General Relativity due to Ehlers, Pirani and Schild
\cite{4}. Moreover, such spaces represent convenient relativistic models
because they have the same conformal and projective properties as the Riemannian
space $(M, \gamma_{ij})$.

Denoting by $r^i_{jkl}$ the curvature tensor field of the metric $\gamma_{ij}$,
by $\gamma^{ij}$ the inverse matrix of $\gamma_{ij}$, $r_{ij}=r^k_{ijk}$,
$r=\gamma^{ij}r_{ij}$,
$\delta/\delta x^i=\partial/\partial x^i-N^j_i(\partial/\partial y^j)$,
$\sigma_i=\delta\sigma /\delta x^i$ and
$\dot\sigma_i=\partial\sigma /\partial y^i$, we shall use the
following notations
\begin{center}
$\sigma^H=\gamma^{kl}\sigma_k\sigma_l$, $\sigma_{ij}=\sigma_{i\vert j}+\sigma_i
\sigma_j-\gamma_{ij}\sigma^H/2$, $\bar\sigma=\gamma^{ij}\sigma_{ij}$,\\
$\sigma^V=\gamma^{ab}\dot\sigma_a\dot\sigma_b$, $\dot\sigma_{ab}=\dot\sigma_
a\vert_b+\dot\sigma_a\dot\sigma_b-\gamma_{ab}\sigma^V/2$, $\dot\sigma=\gamma^
{ab}\dot\sigma_{ab}$,
\end{center}
where ${_\vert}_i$ (resp. $\vert_a$) represents the $h$- (resp. $v$-) covariant
derivative induced by the non-linear connection $N^i_j$.

Developping the formalism presented in \cite{4,5}, the folllowing Maxwell's
equations hold
$$\left\{\begin{array}{lll}
F_{ij\vert k}+F_{jk\vert i}+F_{ki\vert}=\sum_{(ijk)}g_{ip}r^h_{qjk}\dot\sigma
_hy^py^q,\\
F_{ij}\vert_k+F_{jk}\vert_i+F_{ki}\vert_j=-(f_{ij\vert k}+f_{jk\vert i}+
f_{ki\vert}),\\
f_{ij}\vert_k+f_{jk}\vert_i+f_{ki}\vert_j=0,
\end{array}\right.
$$
where the electromagnetic tensors $F_{ij}$ and $f_{ij}$ are
$$
F_{ij}=(g_{ip}\sigma_j-g_{jp}\sigma_i)y^p,\;
f_{ij}=(g_{ip}\dot\sigma_j-g_{jp}\dot\sigma_i)y^p.
$$

Also, the Einstein's equations will take the form
$$
\left\{\begin{array}{ll}\medskip
\displaystyle{r_{ij}-{1\over 2}r\gamma_{ij}+t_{ij}={\cal K}T^H_{ij}},\\
(2-n)(\dot\sigma_{ab}-\dot\sigma\gamma_{ab})={\cal K}T^V_{ab},
\end{array}\right.
$$
where $T^H_{ij}$ and $T^V_{ab}$ are the $h$- and $v$- components of the energy
momentum tensor field, ${\cal K}$ is the gravific constant and
$$
t_{ij}=(n-2)(\gamma_{ij}\bar\sigma-\sigma_{ij})+\gamma_{ij}r_{st}y^s\gamma^{tp}
\dot\sigma_p+\dot\sigma_ir^a_{tja}y^t-\gamma_{is}\gamma^{ap}\dot\sigma_p
r^s_{tja}y^t.
$$

{\bf Remark.} For the form of generalized Einstein-Yang-Mills equations in a
space $(M,e^{2\sigma(x,y)}\gamma_{ij}(x))$, see \cite{1,5}.

Consequently, in certain particular case, it is posible to build a generalized
Lagrange geometry naturally attached to a system of PDEs.

{\bf Open problem.} Is it possible to build a unique generalized Lagrange geometry
naturally asociated to a given PDEs system, in the large?

\begin{tabbing}
\kill
\end{tabbing}

\begin{center}
University POLITEHNICA of Bucharest\\
Department of Mathematics I\\
Splaiul Independentei 313\\
77206 Bucharest, Romania\\
e-mail:udriste@mathem.pub.ro\\
e-mail:mircea@mathem.pub.ro
\end{center}

\end{document}